\documentclass{article}

\usepackage{graphicx}
\usepackage{caption}
\usepackage{subcaption}

\usepackage{amsmath,amsfonts,amssymb}
\usepackage[margin=3.5cm]{geometry}
\usepackage{tikz}
\usetikzlibrary{arrows}
\usepackage{latexsym}
\usepackage{epsfig,overpic}
\usepackage{stmaryrd}
\usepackage{wrapfig}

\let\ep\varepsilon
\newcommand{\R}{\mathbb R}

\newcommand{\bH}{\mathbf H}
\newcommand{\bI}{\mathbf I}

\newcommand{\blf}{\mathbf f}

\newcommand{\bn}{\mathbf n}
\newcommand{\bh}{\mathbf h}
\newcommand{\bp}{\mathbf p}
\newcommand{\bq}{\mathbf q}
\newcommand{\br}{\mathbf r}
\newcommand{\bs}{\mathbf s}

\newcommand{\bx}{\mathbf x}

\newcommand{\T}{\mathcal T}

\newcommand{\uDelta}{{\Delta}_{\Gamma}}

\def\dO{{\partial\Omega} }
\def\div{\operatorname{div} }
\newtheorem{remark}{Remark}

\begin{document}

\title{Numerical integration over implicitly defined domains for higher order unfitted finite element methods\thanks{Partially supported by NSF through the Division of Mathematical Sciences grants 1315993 and 1522252.}}
\author{
Maxim A. Olshanskii\footnotemark[1]
\and
Danil Safin\footnotemark[1]
}
\footnotetext[1]{Department of Mathematics, University of Houston, Houston, Texas 77204-3008 {\tt (molshan,dksafin)@math.uh.edu}}
\date{}

\maketitle

\begin{abstract}
The paper studies several approaches to numerical integration over a domain defined implicitly by an indicator function such as the level set function. The integration methods are based on subdivision, moment--fitting, local quasi-parametrization and Monte-Carlo  techniques. As an application of these techniques, the paper addresses numerical solution of elliptic PDEs posed on domains and manifolds defined implicitly. A higher order unfitted finite element method (FEM) is assumed for the discretization.  In such a method the underlying mesh is not fitted to the geometry, and hence the errors of numerical integration over curvilinear elements affect the accuracy of the finite element solution together with approximation errors. The paper studies the numerical complexity of the integration procedures  and the performance of unfitted FEMs which employ these tools.
\end{abstract}


\section{Introduction} Numerical approaches for solving PDEs that integrate the underlying geometric information, such as isogeometric analysis \cite{hughes2005isogeometric}, are in the focus of research over the last decade. In the isogeometric analysis, the functions used for geometry representation are also employed to define functional spaces
in the Galerkin method. Several other approaches try to make the grid generation and geometry description
independent.  Unfitted finite element methods, such as immersed boundary methods \cite{mittal2005immersed}, extended FEM \cite{belytschko1999elastic,fries2010extended}, cut FEM~\cite{burman2014cutfem}, or trace FEM~\cite{olshanskii2009finite}, use a sufficient regular background grid, but account for geometric details by modifying the spaces of test and trial functions or the right-hand side functional.

The accuracy of unfitted FEM  depends on several factors. These are the approximation property of basic finite element space, the accuracy of underlying geometry recovering, and the error introduced by numerical integration. The present paper first reviews the analysis of two higher order  finite element methods:
one is an unfitted FEM for the Neumann problem in a bounded curvilinear domain,
another one is a narrow-band FEM for an elliptic PDE posed on a closed smooth manifold. In both cases, the domain (a volume or a surface) is given implicitly by a discrete level-set function. Practical implementation of these
methods (as well as many other unfitted FEM) leads to the following problem: Given a simplex $K\in\mathbb{R}^N$, a smooth function $f$ defined on $K$ and a polynomial $\phi_h$  of degree $q$ such that $|\nabla \phi_h|\ge c_0>0$ on $K$, evaluate the integral
\begin{equation}\label{Int}
  I_{K}(\phi_h,f):=\int_{Q}f\,d\bx,\quad \text{with}~Q=\{\bx\in K\,:\,\phi_h(\bx)>0\}.
\end{equation}
The paper focuses on several numerical approaches to problem \eqref{Int} in $\mathbb{R}^2$. In the context of numerical solution of PDEs, $f$ is typically a polynomial;  $q+1$ is the order of geometry recovery.

We note that for $q=1$ one has to integrate $f$ over a polygon (or polyhedron). Then for a polynomial function $f$ an exact numerical integration is straightforward through subdividing the  polygon (polyhedron) into a finite number of triangles (tetrahedra) and applying standard Gauss  quadratures on each triangle. However, for $q>1$ the problem appears to be less trivial and  building an exact quadrature rule for  $I_{K}(\phi_h,f)$ does not look feasible. This can be realised by considering the simple 1D example with $K=(0,1)$ and $f\equiv1$.
Solving \eqref{Int} becomes equivalent to finding the root of  $\phi_h\in (0,1)$. The latter problem is resolved (only) in radicals for $2\le q\le 4$ and is well known to have no general algebraic solution for $q>4$ by the Abel theorem.

The problems of building numerical quadratures for the implicitly defined volume integrals \eqref{Int} and the implicitly defined surface integrals $\int_{S}f\,d\bs,\quad \text{with}~S=\{\bx\in K\,:\,\phi_h(\bx)=0\}$ have been already addressed in the literature and several
techniques have been applied in the context of XFEM and other unfitted FE methods. One straightforward approach consists of employing the smeared Heaviside function $H_\ep$, see, e.g., \cite{QuadHeav}. Then for the regularized problem,  one applies a standard Gaussian quadrature rule on the simplex $K$ with weights $\omega_i$ and nodes $\bx_i$:
\[
I_{K}(\phi_h,f)=\int_{K}f H(\phi_h)\,d\bx\approx \int_{K}f H_\ep(\phi_h)\,d\bx\approx \sum_i \omega_i f(\bx_i) H_\ep(\phi_h(\bx_i)).
\]
However, for a general superposition of $K$ and the zero level set of $\phi_h$,  the smearing leads to significant integration errors which are hard to control. Another numerical integration technique is based on an approximation of $Q$ by elementary shapes. Sub-triangulations
or quadtree (octree) Cartesian meshes are commonly used for these purposes. On each elementary shape a standard quadrature rule is applied.  The sub-triangulation is often adaptively refined towards the zero level of $\phi_h$. The approach is popular in combination with higher order extended FEM for problems with interfaces, see, e.g., \cite{abedian2013performance,moumnassi2011finite,dreau2010studied},
and the level-set method~\cite{min2007geometric,holdych2008quadrature}.
Although numerically stable, the numerical integration based on sub-partitioning may significantly increase the computational complexity of a higher order finite element method, since the number of function evaluations per $K$ scales with $h^{-p}$ for some $p>0$ depending on the order of the FEM. In several recent papers~\cite{muller2013highly,saye2015high,fries2015higher} techniques for numerical integration over implicitly defined domains were devised that have optimal computational complexity. The moment--fitting method from ~\cite{muller2013highly} uses polynomial divergence free basis of vector function to approximate the integrand and further reduce the volume integrals to boundary integrals to find the weights of a quadrature formula by a least-square fitting procedure. We recall the  moment--fitting method in section~\ref{s_MF}.
For the case when $K$ is a hyper-rectangle, the approach in \cite{saye2015high} converts the implicitly given geometry into the graph of an implicitly defined height function. The approach leads to a recursive algorithm on the number of spatial dimensions which requires only one-dimensional root finding and one-dimensional Gaussian quadrature. In \cite{fries2015higher}, a zero level-set of $\phi_h$ is approximated by  higher order  interface elements. These elements are extended inside $K$ so that $K$ is covered by regular and curvilinear simplexes.
For curvilinear simplexes a mapping to the reference simplex is constructed. Further standard Gauss quadratures are applied.

In this paper, we  develop an approach for~\eqref{Int} based on the local quasi-parametrization of the zero level set of $\phi_h$.
Similar to  \cite{saye2015high} the zero level set is treated as a graph of an implicitly given function.  With the help of a 1D root finding procedure the integration over $Q$ is reduced to the integration over regular triangles and the recursive application of 1D Gauss quadratures. The technique is also related to the method of local parametrization for higher order surface finite element method of  \cite{grande2014higher}. We  compare the developed method with several other approaches for the numerical integration in the context
of solving partial differential equations in domains with curvilinear boundaries and over surfaces. For the comparison purpose we consider
the method of sub-triangulation for $Q$, the moment--fitting method, and the Monte-Carlo method.

The rest of the paper is organized as follows. We first recall  unfitted finite element method for solving an elliptic PDE in a domain with curvilinear boundary and an elliptic PDE posed on a surface. For the surface PDE we use the method from~\cite{OlshSafin} of a regular extension to a narrow band around the surface. The error analysis of these unfitted FE methods is also reviewed.
Further in section~\ref{s_int} we discuss the methods for numerical integration of \eqref{Int}, which further used to build the FEM stiffness matrices. Section~\ref{s_num} collects the result of numerical experiments. Section~\ref{s_concl} concludes the paper with a few closing remarks.

\section{Unfitted FEM}\label{s_FEM}

We assume that $\Omega$ is an open bounded subset in $\R^N$, $N=2,3$,  with a boundary  $\Gamma$, which is a connected $C^2$ compact hypersurface in $\R^N$.
In this paper we apply  unfitted FE methods to  elliptic  equations posed in $\Omega$ and on  $\Gamma$. As model problems, let us consider the Poisson and the  Laplace--Beltrami  problems:
\begin{equation}\label{Poiss}
\begin{aligned}
-\Delta u + \alpha\, u &=f\quad\text{in}~\Omega,\\
\frac{\partial u}{\partial n}&=0\quad\text{on}~\Gamma,
\end{aligned}
\end{equation}
and
\begin{equation}\label{LBeq}
-\uDelta u + \alpha\, u =g\quad\text{on}~\Gamma,
\end{equation}
with some strictly positive $\alpha\in L^\infty(\Omega)$ or  $\alpha\in L^\infty(\Gamma)$, respectively.

\subsection{Preliminaries}
To define finite element methods, we need a formulation of the surface PDE \eqref{LBeq} based on normal extension to a narrow band. First, we introduce some preliminaries. Denote by $\Omega_d$ a domain consisting of all points within a distance from $\Gamma$ less than some $d>0$:
\begin{equation*}\label{Omega_d}
\Omega_d = \{\, \bx \in  \R^3~:~{\rm dist}(\bx,\Gamma) < d\, \}.
\end{equation*}
Let $\phi: \Omega_d \rightarrow \R$ be the
signed distance function, $|\phi(x)|:={\rm dist}(\bx,\Gamma)$ for all
$\bx \in \Omega_d$.  The surface $\Gamma$ is the zero level
set of $\phi$:
\begin{equation*}
\Gamma=\{\bx\in\mathbb{R}^3\,:\,\phi(\bx)=0\}.
\end{equation*}
We may assume $\phi < 0$ on the interior of $\Gamma$  and $\phi >0$ on the exterior.
We define $\bn(\bx):=\nabla \phi(\bx)$ for all
$\bx \in \Omega_d$. Thus, $\bn$ is the normal vector  on $\Gamma$,  and $|\bn (\bx)|=1$ for all $\bx\in \Omega_d$.  The Hessian of $\phi$ is denoted by $\bH$:
\begin{equation*}
  \bH(\bx)=\mathrm{D}^2\phi(\bx) \in \R^{3 \times 3} \quad \text{for all} ~~\bx \in \Omega_d.
\end{equation*}
For  $\bx \in \Gamma$, the non-zero eigenvalues of $\bH(\bx)$  are the principal curvatures. Hence, one can choose such sufficiently small positive $d=O(1)$ that
$\bI-\phi\bH$ is uniformly positive definite on $\Omega_d$.
For $\bx\in\Omega_d$ denote by $\bp(\bx)$ the closest point on $\Gamma$.
%
Assume that $d$ is sufficiently small such that the decomposition $\bx=\bp(\bx)+ \phi(\bx)\bn(\bx)$ is unique for all $\bx \in \Omega_d$.
 For a  function $v$ on $\Gamma$ we define its extension to $\Omega_d$:
\begin{equation*} \label{extension}
 v^e(\bx):= v(\bp(\bx)) \quad \text{for all}~~\bx \in \Omega_d.
\end{equation*}
Thus, $v^e$ is the extension of $v$ along normals on $\Gamma$.

We look for $u$ solving the following elliptic problem
 \begin{equation}
\label{ExtNew}
\begin{split}
-\div\mu(\bI-\phi\bH)^{-2}\nabla u+\alpha^e\mu\, u&=g^e\mu\quad \text{in}~~\Omega_d,\\
\frac{\partial u}{\partial \bn}&=0\qquad \text{on}~~\dO_d,
\end{split}
\end{equation}
with $\mu=\text{det}(\bI-\phi\bH)$.
The Neumann condition in \eqref{ExtNew} is the natural boundary condition.  The following  results about the well-posedness of \eqref{ExtNew} and its relation to the surface equations \eqref{LBeq} have been proved in \cite{OlshSafin}:
\begin{description}
\item{}(i) For $g\in L^2(\Gamma)$, the problem~\eqref{ExtNew} has the unique weak solution
$u\in H^1(\Omega_d)$, which satisfies  $\|u\|_{H^1(\Omega_d)}\le C\,\|g^e\|_{L^2(\Omega_d)}$,
with a constant $C$ depending only on $\alpha$ and $\Gamma$;\\[-3ex]
\item{}(ii) For the solution $u$ to \eqref{ExtNew} the trace function $u|_{\Gamma}$ is
an element of $H^1(\Gamma)$ and solves a weak formulation of the surface equation \eqref{LBeq}.\\[-3ex]
\item{}(iii) The solution $u$ to \eqref{ExtNew} satisfies $(\nabla\phi)\cdot (\nabla u)=0$. Using the notion of normal extension, this can be written as
$
u = (u|_\Gamma)^e~ \text{in}~\Omega_d;
$\\[-3ex]
\item{}(iv)  Additionally assume $\Gamma\in C^3$, then $u\in H^{2}(\Omega_d)$ and
$
\|u\|_{H^2(\Omega_d)}\le C\,\|g^e\|_{L^2(\Omega_d)},
$
with a constant $C$ depending only on  $\alpha$, $\Gamma$ and $d$.
\end{description}

\subsection{FEM formulations}\label{s_FEM}
Let  $\Omega^{\rm bulk}\subset\mathbb{R}^N$, $N=2,3$, be a polygonal (polyhedral) domain such that $\Omega\subset\Omega^{\rm bulk}$ for problem $\eqref{Poiss}$
and $\Omega_d\subset\Omega^{\rm bulk}$ for problem $\eqref{LBeq}$.
Assume we are given a family $\{\T_h\}_{h>0}$ of regular triangulations of $\Omega^{\rm bulk}$ such that $\max_{T\in\T_h}\mbox{diam}(T) \le h$.
For a triangle  $T$ denote by $\rho(T)$  the diameter of the inscribed circle.
Denote
\begin{equation}\label{beta}
\beta=\sup_{T\in\T_h}\mbox{diam}(T)/\inf_{T\in\T_h}\rho(T)\,.
\end{equation}
For the sake of presentation, we assume that triangulations of $\Omega^{\rm bulk}$ are quasi-uniform, i.e., $\beta$ is uniformly bounded in $h$.

It is computationally convenient  not to align (not to fit) the mesh  to $\Gamma$ or  $\dO_d$.
Thus, the computational domain $\Omega_h$ approximates $\Omega$ or $\Omega_d$ and has a piecewise smooth boundary which is \textit{not fitted} to the mesh $\T_h$.

Let $\phi_h$ be a continuous piecewise, with respect to $\T_h$, polynomial approximation  of the surface distance function $\phi$ in the following sense:
\begin{equation}\label{phi_h}
\|\phi-\phi_h\|_{L^\infty(\Omega)}+ h \|\nabla(\phi-\phi_h)\|_{L^\infty(\Omega)}\le c\,h^{q+1}
\end{equation}
with some $q\ge1$ (for problem \eqref{Poiss}  a generic level-set function $\phi$ can be considered, not necessary a signed distance function). Then one defines
\begin{equation}\label{Omega_h}
\begin{aligned}
\Omega_h=\{\, \bx \in  \R^3~:~\phi_h(\bx) < 0\, \}&\quad \text{for problem \eqref{Poiss}}\\
\Omega_h=\{\, \bx \in  \R^3~:~|\phi_h(\bx)| < d_h\, \}&\quad \text{for problem \eqref{LBeq}},
\end{aligned}
\end{equation}
with $d_h=O(h)$, $d_h\le d$.
Note, that in some applications the surface $\Gamma$ may not be known explicitly and only a finite element  approximation $\phi_h$ to the distance function $\phi$ is known. Otherwise, one may set $\phi_h:=J_h(\phi)$, where  $J_h$ is a suitable piecewise polynomial interpolation operator. Estimate \eqref{phi_h} is reasonable, if $\phi_h$ is a polynomial of degree $q$ and $\phi\in C^{q+1}(\Omega_d)$. The latter is the case for $C^{q+1}$-smooth $\Gamma$.

The space of all continuous piecewise polynomial functions of a degree $r\ge1$ with respect to   $\T_h$ is our finite element space:
\begin{equation}\label{FEspace}
V_h:=\{v\in C(\T_h)\,:~ v|_T\in P_r(T)\quad\forall\,T\in\T_h\}.
\end{equation}
The finite element method for problem \eqref{Poiss} reads:
Find $u_h\in V_h$ satisfying
\begin{equation}\label{FEmeth}
\int_{\Omega_h}\left[\nabla u_h\cdot\nabla v_h + \alpha^e\,u_h v_h\right]\,d\bx =\int_{\Omega_h} f^ev_h\,d\bx\quad\forall\,v_h\in V_h,
\end{equation}
where $\alpha^e, f^e$ are suitable extensions of $\alpha$ and $f$ to $\Omega^{\rm bulk}$.
For problem \eqref{LBeq}, the finite element method is based on the extended formulation \eqref{ExtNew} and consists of
finding $u_h\in V_h$ that satisfies
\begin{equation}\label{FEmethLB}
\int_{\Omega_h}\left[(\bI-\phi_h\bH_h)^{-2}\nabla u_h\cdot\nabla v_h + \alpha^e\,u_h v_h\right]\,\mu_hd\bx =\int_{\Omega_h} g^ev_h\,\mu_hd\bx\quad\forall\,v_h\in V_h.
\end{equation}
It should be clear  that only those basis functions from $V_h$ contribute to the finite element formulations \eqref{FEmeth}, \eqref{FEmethLB} and are involved in computations that do not vanish everywhere on $\Omega_h$.
Error estimates for the unfitted finite element methods \eqref{FEmeth}, \eqref{FEmethLB} depend on the geometry approximation and the order of finite elements.   Let $\Gamma\in C^{r+1}$  and assume $u\in H^{r+1}(\Omega)$ solves the Neumann problem \eqref{Poiss} and $u_h\in V_h$ solves \eqref{FEmeth}. Then it holds
\begin{equation}\label{FEerror1}
\|u^e-u_h\|_{L^2(\Omega_h)}+h\|u^e-u_h\|_{H^1(\Omega_h)} \le C\,(h^{r+1}+h^{q+1}),
\end{equation}
where a constant $C$ is independent of $h$,   $r$ is the degree of the finite element polynomials, and $q+1$ the geometry approximation order defined in \eqref{phi_h}.

For the Neumann problem with $\alpha=0$ and a compatibility condition on the data, the estimate in \eqref{FEerror1} is proved in \cite{barrett1987practical} subject to an additional assumption on
$\Omega_h$. In  that paper it is assumed that $\partial \Omega_h$ and $\Gamma$ match on the edges of $K\in\T_h$ intersected by $\Gamma$. This is not necessary the case for $\Omega_h$ defined implicitly from the discrete level function $\phi_h$. For implicitly defined domains, the convergence result in \eqref{FEerror1} follows from the more recent analysis in ~\cite{gross2015trace}, where a suitable mapping  of $\Omega_h$ on $\Omega$ was constructed.


For the FE formulation \eqref{FEmethLB} of the Laplace-Beltrami problem \eqref{LBeq}, approximations to $\phi$ and $\bH$ are required.
If $\Gamma$ is given explicitly, one can compute $\phi$ and $\bH$ and set  $\phi_h=\phi$, $\bH_h=\bH$ and
$\mu_h=\mbox{det}(\bI-\phi_h\bH_h)$ in \eqref{FEmeth}.
Otherwise, if  the surface $\Gamma$ is known approximately as, for example, the zero level set of a finite element distance function $\phi_h$, then, in general, $\phi_h\neq \phi$  and one has to define a discrete Hessian $\bH_h\approx \bH$ and also set $\mu_h=\mbox{det}(\bI-\phi_h\bH_h)$. A discrete Hessian $\bH_h$  can be obtained from $\phi_h$ by a  recovery method, see, e.g.,~\cite{Hessian1,Hessian0}.
Assume that some $\bH_h$ is provided and denote by $p\ge 0$ the approximation order for $\bH_h$ in the (scaled) $L^2$-norm:
\begin{equation} \label{Hh}
|\Omega_h|^{-\frac12}\|\bH-\bH_h\|_{L^2(\Omega_h)}\le c h^p,
\end{equation}
where $|\Omega_h|$ denotes the area (volume) of $\Omega_h$.

The convergence of the finite element method \eqref{FEmeth} is summarized in the following result from~\cite{OlshSafin}.
Let $\Gamma\in C^{r+2}$, $f\in L^\infty(\Gamma)$, and assume $u\in W^{1,\infty}(\Gamma)\cap H^{r+1}(\Gamma)$ solves the surface problems \eqref{LBeq} and $u_h\in V_h$ solves \eqref{FEmeth}. Then it holds
\begin{equation}\label{FEerror2}
\|u-u_h\|_{H^1(\Gamma)} \le C\,(h^r+h^{p+1}+h^q),
\end{equation}
where a constant $C$ is independent of $h$, and $r\ge1$, $p\ge0$, $q\ge1$ are the finite elements, Hessian recovery, and distance function approximation orders defined in \eqref{FEspace}, \eqref{phi_h} and \eqref{Hh}, respectively. Numerical experiments in \cite{OlshSafin} show that
$\|u-u_h\|_{L^2(\Gamma)}$ typically demonstrates a one order higher convergence rate than the $H^1(\Gamma)$ norm of the error.

Note that both error estimates \eqref{FEerror1} and \eqref{FEerror2} assume exact numerical integration. In the introduction, we discussed that for $q>1$ the exact numerical integration is not feasible. The error of the numerical quadrature should be consistent with
the finite element interpolation and geometric error to ensure that the FEM preserves the optimal accuracy.

\begin{remark}\rm
For the purpose of improving algebraic properties of an unfitted finite element method a stabilization procedure was suggested in
\cite{burman2012fictitious} for elliptic equations posed in volumetric domains and further extended to surface PDEs in  \cite{Alg1}.
The procedure consists of adding a special term penalizing the jump of the solution gradient over the edges (faces) of the triangles (tetrahedra) cut by $\dO_h$.

Let $\T^h_\Gamma:=\{K\in \T_h\,:\,\text{meas}_{N-1}(K\cap\dO_h)>0\}$ ($\T^h_\Gamma$ is the set of all elements having non-empty intersection with the boundary of the numerical domain $\Omega_h$). By $\mathcal{F}_\Gamma^h$ denote the set of all edges in 2D or faces in 3D shared by any two elements from  $\T^h_\Gamma$. Define the term
\[
J(u,v)=\sum_{F\in\mathcal{F}_\Gamma^h}\sigma_F \int_{F}\llbracket\bn_F\cdot\nabla u\rrbracket\llbracket\bn_F\cdot\nabla v\rrbracket.
\]
Here $\llbracket\bn_F\cdot\nabla u\rrbracket$ denotes the jump of the normal derivative of $u$ across $F$;
$\sigma_E\ge0$ are stabilization parameters (in our computations we set $\sigma_F=1$).  
The edge-stabilized trace finite element reads: Find $u_h\in V_h$ such that
\begin{equation}
a_h(u_h,v_h) + J(u_h,v_h) =f_h(v_h),  \label{FEM_edge}
\end{equation}
for all $v_h\in V_h$. Here $a_h(u_h,v_h)$ and $f_h(v_h)$ are the bilinear forms and the right-hand side functional corresponding to
the finite element methods \eqref{FEmeth} or  \eqref{FEmethLB}.

For $P_1$ continuous bulk finite element methods on quasi-uniform regular tetrahedral meshes,  the optimal orders of convergence for \eqref{FEM_edge} were proved in \cite{burman2012fictitious}.
In our numerical studies we tested the stabilized formulation \eqref{FEM_edge} with higher order elements.
We observed very similar convergence results for the formulations with and without  $J(u_h,v_h)$ term, including sub-optimal/irregular behaviour with moment-fitting quadratures and optimal with other integration techniques. For the systems of linear  algebraic equations we use exact `backslash' solves in either case. Therefore, we shall not report results for \eqref{FEM_edge} in addition to finite element formulations \eqref{FEmeth} and \eqref{FEmethLB}.
\end{remark}

\section{Numerical integration} \label{s_int} \label{s_MF}

For a simplex $K\subset\mathbb{R}^N$, $K\in \T_h$ we are interested in computing integrals  over $Q=\{\bx\in K\,:\,\phi_h(\bx)>0\}$ with a certain accuracy $O(h^{m})$, i.e. for a sufficiently smooth $f$ we look for a numerical quadrature $I_{h,K}(\phi_h,f)$ such that
\begin{equation}\label{aux1}
|I_{K}(\phi_h,f)-I_{h,K}(\phi_h,f)| \le c\,h^{m},\quad m=\min\{q,r\}+N,
\end{equation}
with a constant $c$ uniformly bounded over $K\in\T_h$. We restrict ourselves with the two-dimensional case, $N=2$.

As a pre-processing step we may compute a simple polygonal approximation $Q_K$ to the curvilinear integration domain $Q$. To find the   polygonal subdomain $Q_K$, we invoke a root finding procedure (several iterations of the secant method are used in our implementation) to find all  intersection points  of the zero level set of  $\phi_h(x)$ with the edges of $K$. $Q_K$ is defined as a convex hull of these intersection points and the vertices of $K$ lying in $\Omega_h$.  The integration of a polynomial function $f$ over $Q_K$ can be done exactly using Gaussian quadratures
on a (macro) triangulation of $Q_K$. The integration over $\widetilde{Q}=Q\triangle Q_K$ has to be done approximately.  
We can write
\[
\int_{Q}f\,d\bx=\int_{Q_K}f\,d\bx+\int_{\widetilde{Q}}\text{sign}(\phi_h)f\,d\bx.
\]
For the numerical integration over $\widetilde{Q}$ and $Q$, we consider several approaches.  We start with  the most straightforward technique, the  Monte-Carlo method.
\medskip

%
%

For the  \textit{Monte-Carlo method} we seed $M\gg1$  points $\{\bx_i\}$ using a uniform random distribution over a narrow rectangular strip $S$ containing $\widetilde{Q}$. Further compute
\begin{equation}\label{Mont}
\int_{\widetilde{Q}}\text{sign}(\phi_h) f\,d\bx \approx \frac{|S|}{M}\sum_{i=1}^M g(\bx_i),\quad g(\bx)=f(\bx)\text{sign}(\phi_h(\bx))\chi_{\widetilde{Q}}(\bx).
\end{equation}
Further we calculate the variance $|S|\sqrt{\left|(\sum_{i=1}^M g(\bx_i))^2/M - \sum_{i=1}^M(g(\bx_i)^2)\right|}/M$. If the variance exceeds a predefined   threshold $\ep=O(h^m)$, then we increase the set of points used and update the sum on the right-hand side of \eqref{Mont}. Noting $|S|=O(h^3)$,  a conservative estimate
gives $M\simeq O(h^{3-2m})$.

In the Monte-Carlo method  the number of function evaluations per a grid element intersected by $\dO_h$  is too far from
being optimal. A {sub-triangulation} method below allows to decrease the number of function evaluations per cell and still constitutes a very robust approach.
\medskip

\begin{figure}
\centering
\includegraphics[width=0.45\textwidth]{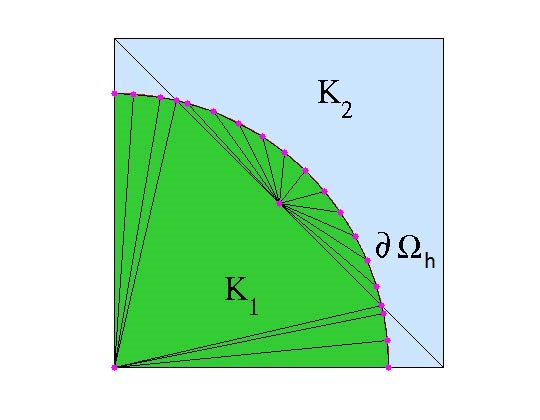}\qquad
\includegraphics[width=0.45\textwidth]{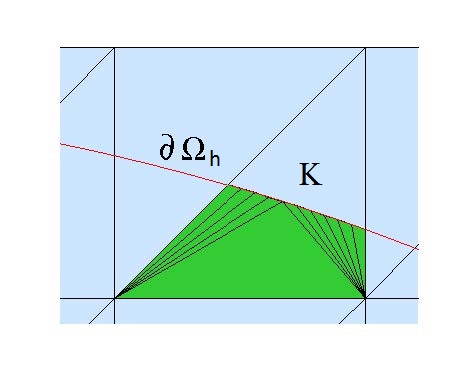}
\caption{Local subdivisions of cut triangles.\label{fig:P2}}
\end{figure}

In the \textit{sub-triangulation} algorithm we construct a local triangulation of ${Q}$. This is done by finding $O(h^{3-m})$ points on the curvilinear boundary  $\dO_h$. In our implementation, these points are found as intersections of a uniform ray corn tailored to a basis point on $\partial Q_K$, see Figure~\ref{fig:P2} for the example of how such local triangulations were constructed for  cut triangles $K_1$ and $K_2$ (left plot) and a cut triangle $K$ (right plot) of a bulk triangulation (FE functions are integrated over the green area).  Further, the integral over a cut element is computed as the
sum of integrals over the resulting set of smaller triangles.
The approach can be viewed as building a local piecewise linear approximation of $\phi_h$ with some $h'=O(h^{m'})$,
with $m'=\min\{q,r\}$, i.e. $h'=O(h^2)$ for $P_2$ elements and  $h'=O(h^3)$ for $P_3$ elements.

The number of function evaluations per  a triangle intersected by $\dO_h$ is O($h^{3-m}$), which is better than with the Monte-Carlo method, but still sub-optimal. The first integration method delivering optimal complexity we consider is the  Moment-Fitting method from~\cite{muller2013highly}.
\medskip

\begin{wrapfigure}{l}{0.5\textwidth}
\centering
\includegraphics[width=0.5\textwidth]{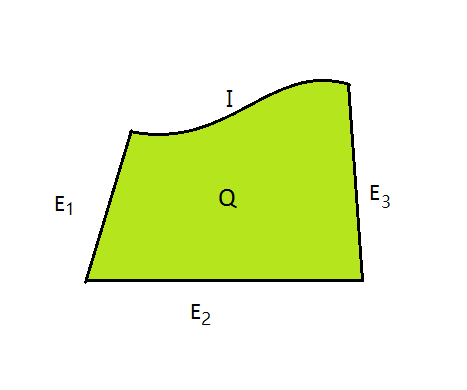}
\caption{Integration domain $Q$ for the moment--fitting.\label{fig:MF}}
\end{wrapfigure}

In the \textit{Moment-Fitting method}, one first defines a set of points $\{\bx_i\}$, $i=1,\dots,M$,  for a given cell $K$ intersected by $\dO_h$. The choice of the points can be done for a reference triangle and is independent on how $\dO_h$ intersects $K$:  $\{\bx_i\}$ can be regularly spaced, come from a conventional quadrature scheme, or even randomly distributed.   The boundary of the integration domain $Q$ consists of  straight edges $E_k$ and the curvilinear part $I$, cf. Figure~\ref{fig:MF}. Let $\bn_k$, $\bn_I = (\nabla \phi_h)/|\nabla \phi_h|$ to be the unit outward normal vectors for each $E_k$ and $I$ respectively.
The moment--fitting method calculates quadrature weights for the chosen points $\{\bx_i\}$ by finding the least-square solution to the system
\begin{equation}\label{LSM}
\sum_{i=1}^{M} \omega_i g_j(\bx_i)=\int_Q g_j\,d\,\bx
\end{equation}
for a given set of basis functions $ \{g_j\}$, $j=1,\dots,K$.
For example, for P2 FEM we use the basis
\begin{equation}\label{func_g} \mathcal{G} = \{g_j\} = \{ 1, \ x, \ y, \ x^2, \ xy, \ y^2 \}
\end{equation}
The number $M$ of points $\{\bx_i\}$ is recommended in~\cite{muller2013highly}  to exceed the number of basis functions, i.e.,  $K< M$ holds.

The integrals on the right-hand side of \eqref{LSM} are evaluated approximately by resorting to divergence-free basis functions.
For divergence-free basis functions, one applies the divergence theorem to reduce dimensions of the integrals and relate integration over implicit surface to simple line integrals. In 2D, the  div-free basis complementing \eqref{func_g} is given by
$$
\mathcal{F} = \{\blf_j\} = \left\{ \begin{array}{lllllllll}
1 & 0 & 0 & x & y & y^2 & 2xy & x^2 & 0 \\
0, & 1, & x, & -y, & 0, & 0, & -y^2, & -2xy, & x^2 \\
\end{array}
\right\}
$$
Since $\partial Q$ is a closed piecewise smooth curve, for $\blf_j\in\mathcal{F}$ it holds
\[ 0 = \int_Q \mathrm{div}\, \blf_j(\bx) d\,\bx=
\int_{\partial Q} \blf_j\cdot \bn_{\partial Q} ds = \int_{I} \blf_j\cdot \frac{\nabla \phi_h}{|\nabla \phi_h|} ds + \sum_k \int_{E_k} \blf_j\cdot \bn_k ds
\]
Hence, we obtain
\[
\int_{I} \blf_j\cdot \frac{\nabla \phi_h}{|\nabla \phi_h|} ds  = - \sum_k \int_{E_k} \blf_j\cdot \bn_k ds
\]
Integrals over any interval $E_k$ on the right side can be computed with a higher order Gauss quadrature rule. This allows to build a quadrature for the numerical integration over implicitly given curvilinear edge $I$ based on the interior nodes  $\{\bx_i\}$.
To this end, one calculates  the weights $\{v_i\}$ for the set of nodes $\{\bx_i\}$ by solving the following  system:
\[
 \sum_{i=1}^{M}  \blf_j(\bx_i)\cdot \frac{\nabla \phi_h(\bx_i)}{|\nabla \phi_h(\bx_i)|} v_i = \int_I \blf_j\cdot\frac{\nabla \phi_h}{|\nabla \phi_h|}  ds,\quad j=1,\dots,K.
\]
For area integration, take the second set of functions related to $\mathcal{G}$ as $ \mathrm{div}\,\bh_j = 2g_j$:
$$ \mathcal{H} = \{\bh_j\} = \left\{ \begin{array}{ccccccc}
x & x^2/2 & xy  & x^3/3 & x^2y/2 & xy^2 \\
y & xy & y^2/2  & x^2y & xy^2/2 & y^3/3 \\
\end{array}
\right\} $$
The divergence theorem gives
$$
2\int_Q g_j d\bx = \int_{\partial Q} \bh_j\cdot \bn_{\partial Q} ds = \int_{I} \bh_j\cdot \frac{\nabla \phi_h}{|\nabla \phi_h|} ds + \sum_k \int_{E_k} \bh_j\cdot \bn_k ds
$$
Now the right-hand side values in \eqref{LSM} are (approximately) computed using the above identity and the previously computed surface quadrature rule for the numerical integration over $I$:
$$
2\int_Q g_j d\bx \approx  \sum_{i=1}^{M}  \bh_j(\bx_i)\cdot \frac{\nabla \phi_h(\bx_i)}{|\nabla \phi_h(\bx_i)|} v_i + \sum_k \int_{E_k} \bh_j\cdot \bn_k ds
$$
 This algorithm can be extended to higher order quadratures by expanding the function sets $\mathcal{G}$, $\mathcal{F}$ and $\mathcal{H}$; or to integration over a higher dimensional curvilinear simplex by adding a further moment--fitting step.

The complexity of the moment--fitting integration is optimal, e.g. O(1) of function evaluations per triangle.
However, the weights computed by the fitting procedure are not necessarily all non-negative. There is no formal prove of
the resulting quadrature accuracy as in \eqref{aux1} or the consistency order.  Furthermore, in experiments we observe
that the finite element methods with stiffness matrices assembled using
 moment--fitting can be less stable compared to applying other numerical integration techniques discussed here.
Below we consider another integration method of optimal complexity.
\medskip

\begin{figure}
\begin{center}
\begin{minipage}{.33\textwidth}
\begin{tikzpicture}[scale=4,cap=round]
  \tikzstyle{axes}=[]
  \tikzstyle{important line}=[very thick]
  \tikzstyle{information text}=[rounded corners,fill=red!10,inner sep=1ex]

  \draw[style=important line] (1, 0) -- (1, 1);

  \draw[style=important line] (0,0) -- (1,0);

  \draw[style=important line] (1,1) -- (0, 0);

  \draw[style=dashed, color=blue] (0.60,0.595) -- (0.8,0);

  \draw (0.60,0.595) circle (.03cm) node[blue, left=5pt,fill=none] {$\bp_1$};
  \draw (0.8,0) circle (.03cm) node[blue, below left=5pt,fill=none] {$\bp_2$} ;

  \draw [red,thick,domain=-10:50] plot ({cos(\x)-0.2}, {sin(\x)}) node[above=2pt,fill=white] {$\phi(x)=c_0$} ;

\end{tikzpicture}
\end{minipage}
\begin{minipage}{.33\textwidth}
\begin{tikzpicture}[scale=5,cap=round]
\tikzset{cross/.style={cross out, draw,
         minimum size=2*(#1-\pgflinewidth),
         inner sep=0pt, outer sep=0pt}}
  \tikzstyle{axes}=[]
  \tikzstyle{important line}=[very thick]
  \tikzstyle{information text}=[rounded corners,fill=red!10,inner sep=1ex]

  \draw[style=solid, color=blue] (0.43,1.25) -- (0.95,0);

  \draw  (0.43,1.25) circle (.02cm) node[blue, left=3pt,fill=white] {$\bp_1$};
  \draw (0.95,0) circle (.02cm) node[blue, below left=3pt,fill=white] {$\bp_2$} ;

  \draw [red,thick,domain=-20:68] plot ({cos(\x)}, {sin(\x)+sin(20)}) ;

  \draw [thick,->] (0.57, 0.9) -- (1, 1.1) node[blue, above right=2pt,fill=white] {};

    \draw [thick,->] (0.83, 0.3) -- (1.26, 0.5) node[blue, above right=2pt,fill=white] {};

  \draw (0.57, 0.9) -- (0.62, 0.92) -- (0.64, 0.88) -- (0.59, 0.86) ;
  \draw (0.83, 0.30) -- (0.88, 0.32) -- (0.90, 0.28) -- (0.85, 0.26) ;

  \draw [fill=green] (0.57, 0.9) circle (.02cm) node[blue, left=5pt,fill=white] {$\bq_1$};
  \draw [fill=green] (0.82, 0.30)  circle (.02cm) node[blue, left=5pt,fill=white] {$\bq_2$};

  \draw [fill=green] (0.77, 0.99)  circle (.02cm) node[blue, above=5pt,fill=none] {$\hat   \bq_1$};
  \draw [fill=green] (1, 0.38)  circle (.02cm) node[blue, below right=5pt,fill=white] {$\hat \bq_2$};

    \draw  (0.43,1.25) circle (.02cm) node[blue, left=3pt,fill=white] {$\bp_1$};

     \draw  (0.91,0.6) circle (.00cm) node[black, left=3pt,fill=white] {$\widetilde{Q}$};
     \draw  (1.05,1.02) circle (.00cm) node[black, left=3pt,fill=white] {$\bn$};


\end{tikzpicture}
\end{minipage}
%
%
%
%
%
%
%
%
%
%
%
\end{center}
\end{figure}

Consider the curvilinear remainder $\widetilde{Q}$ and denote by  $\bp_1$, $\bp_2$ the intersection points of $\dO_h$ with $Q_K$.
If there are more than 2 such points, then the calculations below should be repeated for each of the simply-connected component of $\widetilde{Q}$. Choose  points $\{\bq_i\}$ on $(\bp_1,\bp_2)$ as a nodes of a Gaussian quadrature with weights $\{\omega_i\}$. Consider the normal vector $\bn$ for the  line passing through $\bp_1$ and $\bp_2$. For each point $\bq_i$, one finds the point $\hat \bq_i$ on the boundary such that $\hat \bq_i = \bq_i+\alpha \bn$, $\alpha\in\mathbb{R}$, and $\phi_h(\hat \bq_i) - {c_0} = 0$, where $c_0$ is the $\phi_h$-level value for $\dO_h$, i.e. $c_0=0$ for the FEM \eqref{FEmeth} and $c_0=\pm d_h$ for the FEM \eqref{FEmethLB}. Secant method finds $\hat \bq_i$ up to machine precision within a few steps. Further we employ the same 1D quadrature rule to place points \{$\br_{ij}$\} on each segment ($\bq_i$, $\hat \bq_i$).
The integral over $\widetilde{Q}$ is computed through
\begin{equation}\label{num_int}
\int_{\widetilde{Q}} \text{sign}(\phi_h) f(x) dx\approx I_{h,\widetilde{Q}}(f) = |\bp_1-\bp_2|\sum_i|\bq_i-\hat \bq_i|\text{sign}(\phi_h(\bq_i)) \omega_i \sum_j f(\br_{ij}) \omega_j.
\end{equation}
The weights for integrating over $Q$ may be combined in $\omega_{ij} = |\bp_1-\bp_2| |\bq_i-\hat \bq_i| \text{sign}(\phi_h(\bq_i)) \omega_i \omega_j $ to write the final quadrature formula  $\sum_i \sum_j f(\br_{ij}) \omega_{ij}$. We remark that for problem \eqref{FEmethLB}, the factor $\text{sign}(\phi_h(\bq_i))$ in \eqref{num_int} is replaced by $\text{sign}(\phi_h(\bq_i)+d_h)$ or $\text{sign}(d_h-\phi_h(\bq_i))$ depending on the level set of $\dO_h$.

Similar to the moment--fitting method, the complexity of the numerical integration is optimal, i.e. O(1) of function evaluations per triangle. All weights $\omega_{ij}$ are positive, and the accuracy analysis of numerical integration of a smooth function over $\Omega$, which invokes the constructed quadrature for handling boundary terms, is straightforward and outlined below.
\medskip

%

We estimate the error of integration of a sufficiently smooth $f$ over $\Omega$. The integration uses a conventional quadrature scheme for interior cells of $T_h$ and the quadrature \eqref{num_int} for the curvilinear remainders of cut  cells. This composed numerical integral is denoted by $I_{h,K}(f)$.  Assume the 2D quadrature over regular triangles  has $O(h^m)$ accuracy.
By the triangle inequality, we have
\[
\left\lvert \int_\Omega f \,d\bx - I_{h,\Omega}(f) \right\rvert \leq \left\lvert\int_\Omega f \,d\bx - \int_{\Omega_h} f \,d\bx\right\rvert  + \left\lvert\int_{\Omega_h} f\,d \bx - I_{h,\Omega}(f)\right\rvert
\]
We apply the co-area formula to estimate the first term
\begin{multline*}
\left\lvert\int_\Omega f \,d\bx - \int_{\Omega_h} f \,d\bx\right\rvert
\leq\int\limits_{(\Omega \triangle\Omega_h) } |f| d\bx
\leq \|f\|_{L^\infty} \int\limits_{(\Omega \triangle\Omega_h) } 1 d\bx \\
\leq \|f\|_{L^\infty} \int\limits_{-\|\phi_h-\phi\|_{L^\infty}}^{+\|\phi_h-\phi\|_{L^\infty}} \int\limits_{\{\phi=t\} } |\nabla \phi| d\bx dt
\le Ch^q \|f\|_{L^\infty}.
\end{multline*}
For the second term we estimate
\[
\left\lvert\int_{\Omega_h} f \,d\bx - I_{h,\Omega}(f) \right\rvert  \leq \sum_{K \in T_h} \left\lvert \int_K f \,d\bx - I_{h,K}(f)\right\rvert,
\]
where $I_{h,K}(f)$ is a quadrature we use to integrate $f$ over $K\cap\Omega_h$.
Decomposing the mesh into interior cells $T_h^{int}$ and cut cells $T_h^{\Gamma}$ (those intersected by $\partial\Omega_h$)
and applying triangle inequalities, we obtain:
\[
\left\lvert\int_{\Omega_h} f d\bx - I_{h,\Omega}(f) \right\rvert  \leq  \sum_{K \in T_h^{int} }\left\lvert\int_K f d\bx - I_{h,K}(f)\right\rvert  + \sum_{K \in T_h^{\Gamma} }\left\lvert\int_{K \cap \Omega_h} f d\bx - I_{h,K}(\phi_h, f)\right\rvert.
\]
The interior cells are integrated with the error $Ch^m$ by a conventional method:
\[
\sum_{K \in T_h^{int} }\left\lvert\int_K f dx - I_{h,K}(f)\right\rvert\le Ch^m.
\]
For $Q=K\cap\Omega_h=Q_K\cup\widetilde {Q}_K$, $\widetilde{Q}_K = Q\triangle Q_K$ and polygonal $Q_K$ defined earlier,  we have
\[
 \sum_{K \in T_h^{\Gamma} }\left\lvert\int_{K \cap \Omega_h} f d\bx - I_{h,K}(\phi_h, f)\right\rvert \leq Ch^m + \sum_{K \in T_h^{\Gamma} }|\int_{\widetilde{Q}_K}\text{sign}(\phi_h) f d\bx- I_{\widetilde{Q}_K}(f)|.
\]

The estimate below assumes that a Gaussian 1D quadrature with $P$ nodes is used in the construction of \eqref{num_int}.   For each interval
$(\bp_1,\bp_2)$ we introduce the local orthogonal coordinate system $(s,t)=\bx(s)+t\bn$, where $s:\,(0,|\bp_2-\bp_1|)\to (\bp_1,\bp_2)$ parameterizes the interval. The graph of the zero level of $\phi_h$ is the implicit function $\gamma(s)$ given by $\phi_h(\bx(s)+\gamma(s)\bn)=0$.
We note the identity
\[
\frac{d^{2P}}{ds^{2P}}\int_0^{\gamma(s)}f(s,t)\,dt=\sum_{k=0}^{2P-1}C^{2P}_k \frac{\partial^k f}{\partial s^k} \frac{d^{2P-k}\gamma }{ds^{2P-k}} +\int_0^{\gamma(s)} \frac{\partial^{2P} f}{\partial s^{2P}}(s,t)\,dt.
\]
We assume $f$ and $\gamma(s)$ to be smooth enough that
\begin{equation}\label{dir_est}
 |\frac{d^{2P}}{ds^{2P}}\int_0^{\gamma(s)}f(s,t)\,dt|\le C_f,
\end{equation}
with a constant $C_f$ uniform over all $K\in\T_h^{\Gamma}$ and independent of $h$.
Applying standard estimates for the Gaussian quadratures, we get:
\begin{equation}
\begin{aligned}
|\int_{\widetilde{Q}_K}\text{sign}(\phi_h)f d\bx &- \sum_i \sum_j \omega_{ij} f(\br_{ij})| = | \int_{\bp_1}^{\bp_2} \int_0^{\gamma(s)}f(s,t)dtds - \sum_i \sum_j \omega_{ij} f(\br_{ij})|\\
 &\le |\sum_i \int_{0}^{\gamma(\bq_i)} f\,dt -  \sum_j \omega_{ij} f(\br_{ij})| + C\,|\bp_1-\bp_2| h^{(2P)} \\
 &\le \sum_i C\,|\gamma(\bq_i)| \|f\|_{W^{2P,\infty}} h^{(2P)} + C\, h^{(2P+1)} \le C\,h^{(2P+1)}.\\
\end{aligned}
\end{equation}
For a 2D boundary, the number of unfitted regions $Q$ grows at an order of O($h^{-1}$). Then the entire error for all unfitted regions will be $O(h^{2P})$. 

To satisfy \eqref{aux1}, it is sufficient to set $2P+1=m\le q+2$. Therefore, for sufficiently smooth $f$ the assumption \eqref{dir_est} is valid if $\gamma(s)\in W^{q+1,\infty}$ with the
$W^{q+1,\infty}$-norm uniformly bounded  over all cut triangles and independent of $h$. The latter follows from our assumptions on $\phi_h$ and $\phi$.  Indeed, let $J_h(\phi)$ be a suitable polynomial interpolant for $\phi$ on a cut triangle $K$.
By triangle inequality we have
\begin{equation}\label{phi_est}
\|\phi_h\|_{W^{q+1,\infty}(K)}\le \|\phi_h-J_h(\phi)\|_{W^{q+1,\infty}(K)}+\|\phi-J_h(\phi)\|_{W^{q+1,\infty}(K)}+\|\phi\|_{W^{q+1,\infty}(K)}.
\end{equation}
Applying the finite element inverse inequality, condition \eqref{beta}, estimate \eqref{phi_h} and approximation properties of polynomials we get  for the first term on the right-hand side of
\eqref{phi_est}:
\begin{multline*}
\|\phi_h-J_h(\phi)\|_{W^{q+1,\infty}(K)}\le C\,h^{-q}\|\phi_h-J_h(\phi)\|_{W^{1,\infty}(K)}\\ \le C\,h^{-q}(\|\phi_h-\phi\|_{W^{1,\infty}(K)}+\|\phi-J_h(\phi)\|_{W^{1,\infty}(K)})\le
C\,,
\end{multline*}
with a constant $C$ independent of $K$ and $h$. Estimating the second and the third terms on the right-hand side of
\eqref{phi_est} in an obvious way,  we obtain
\begin{equation}\label{phi_est2}
\|\phi_h\|_{W^{q+1,\infty}(K)}\le C\,,
\end{equation}
with a constant $C$ independent of $K$ and $h$. The desired estimate on the $W^{q+1,\infty}$-norm of $\gamma$ follows from \eqref{phi_est2}, the properties of implicit function and assumptions on $\phi$.

\section{Numerical examples}\label{s_num}

In this section we demonstrate the results of a few  experiments using different numerical integration approaches described in the previous section.

\subsection{Integral of a smooth function}
We first experiment with computing  integral of a smooth function over an implicitly defined domain in $\mathbb{R}^2$.
For the domain we choose the annular region defined by the level set function
$$ \phi(\bx) =  | \,|\bx|-1|-0.1,
\qquad  \Omega=\{\bx\in \mathbb{R}^2 \,:\,\phi(\bx)<0\}.
 $$
For $f$ given in polar coordinates by
\[
 f(r, \theta) =  10^5 \sin(21 \theta) \sin(5 \pi r).
\]
 the exact value $\int_\Omega f\,d\bx$  is known and can be used to test the accuracy of different approaches.

Due to the symmetry, the computational domain is taken to be the square $\Omega^{\rm bulk}=(0,1)^2$. Further,  uniform triangulation with meshes of sizes $h=0.1\times2^{-i}$, $i=0,\dots,8$, are built to triangulate $\Omega^{\rm bulk}$.
To avoid extra geometric error and assess the accuracy of numerical integration, in these  experiments we set $\phi_h=\phi$. All four methods are set up to deliver the local error estimate \eqref{aux1} with $m=4$  or $m=5$. This should lead to $O(h^3)$ and $O(h^4)$ global accuracy, respectively.

Tables~\ref{table1}--\ref{table2.1} demonstrate that all methods demonstrate  convergent results with local parametrization and sub-triangulation being
somewhat more accurate in terms of absolute error values. At the same time, only  moment--fitting  and local parametrization approaches are optimal
in terms of the computational complexity. Here we measure complexity in terms of the number of function evaluations. The total number of function evaluations to compute integrals over cut elements and interior elements is shown.
Monte-Carlo method appears to be the most computationally expensive. Both Monte-Carlo and sub-triangulation methods become prohibitively expensive for fine meshes so that we  make only 4 refining steps with those methods for $m=4$ and only 2 refining steps for $m=5$. The actual CPU timings (not shown) depend on particular implementation. For a Matlab  code we used, the moment--fitting was the fastest among the four tested  for a given $h$.

\begin{table}[ht]\small\center
\begin{tabular}{l|rr|rr|rr|rr}
    h      & MF       & rate & MC       & rate & ST      &rate&  LP      &rate \\ \hline
    0.1000 & 7.41e+01 &      & 5.50e-01 &      & 1.02e-01&    &  1.66e+00&
 \\ 0.0500 & 6.22e+00 & 3.57 & 4.78e-02 &3.52  & 8.53e-03&3.58&  7.06e-03&7.88
 \\ 0.0250 & 9.14e-02 & 6.09 & 1.56e-03 &4.94  & 4.80e-04&4.15&  4.46e-03&0.66
 \\ 0.0125 & 7.96e-03 & 3.52 & 7.93e-04 &0.98  & 1.35e-05&5.15&  5.91e-05&6.24
 \\ 0.0062 & 1.28e-03 & 2.64 & 1.22e-04 &2.70  & 1.17e-06&3.53&  1.18e-05&2.32
 \\ 0.0031 & 9.05e-06 & 7.14 &         &       &         &    &  2.62e-07&5.49
 \\ 0.0016 & 1.11e-05 &-0.29 &       &         &        &    &   3.99e-08&2.72
 \\ 0.0008 & 4.17e-08 & 8.06 &       &         &        &    &   1.91e-09&4.38
 \\ 0.0004 & 1.68e-07 &-2.01 &      &          &         &   &   1.44e-10&3.73
 \\ \end{tabular}
 \caption{The global error and the error reduction rates  for the numerical integration using moment--fitting (MF), Monte-Carlo (MC), sub-triangulation (ST), and local parametrization (LP) algorithms to handle cut elements.
    The table shows results with $m=4$.} \label{table1}
 \end{table}

 \begin{table}[ht]\small \center
\begin{tabular}{l|r|r|r|r}
    h     &MF   &MC     & ST    &LP  \\ \hline
    0.1000&1\,352    &25\,802      & 4\,136      & 423
 \\ 0.0500&3\,718    &93\,491      & 16\,966     & 1\,608
 \\ 0.0250&11\,726   &346\,005     & 68\,942     & 12\,750
 \\ 0.0125&40\,144   &1\,342\,357  & 277\,456    & 42\,192
 \\ 0.0062&146\,926  &40\,645\,578 & 1\,113\,070 & 151\,022
 \\ 0.0031&561\,912  & -        &     & 570\,104
 \\ 0.0016&2\,194\,452 & -      &     & 2\,210\,836
 \\ 0.0008&8\,671\,104 & -      &     & 8\,703\,872
 \\ 0.0004&34\,471\,164 & -     &     & 34\,536\,700
 \\ \end{tabular}
 \caption{The  number of function evaluations  for the numerical integration using moment--fitting (MF), Monte-Carlo (MC), sub-triangulation (ST), and local parametrization (LP) algorithms to handle cut elements.
    The table shows results with $m=4$.} \label{table1.1}
 \end{table}

 \begin{table}[ht]\small\center
\begin{tabular}{l|rr|rr|rr|rr}
h & MF & rate & MC & rate & ST &rate & LP &rate  \\ \hline
    0.1000 & 2.80e+01 &     & 6.49e-02&      & 5.51e-02 &      & 4.16e-02 &
 \\ 0.0500 & 2.79e-01 & 6.65& 8.64e-04& 6.23 & 7.60e-03 &2.86  & 1.13e-04 & 8.52
 \\ 0.0250 & 4.44e-03 & 5.97&         &      & 2.25e-05 &8.40  & 7.18e-06 & 3.98
 \\ 0.0125 & 3.61e-04 & 3.62&         &      &          &      & 9.75e-09 & 9.52
 \\ 0.0062 & 2.11e-05 & 4.10&         &      &          &      & 1.53e-09 & 2.67
 \\ 0.0031 & 9.59e-07 & 4.46&         &      &          &      & 7.73e-12 & 7.63
 \\
 \end{tabular}
\caption{The global error and the error reduction rates for numerical integration using moment--fitting (MF), Monte-Carlo (MC), sub-triangulation (ST), and local parametrization (LP) algorithms to handle cut elements.
    The table shows results with $m=5$.} \label{table2}
 \end{table}

\begin{table}[ht]\small\center
\begin{tabular}{l|rr|rr|rr|rr}
h  &MF   &MC     & ST    &LP \\ \hline
    0.1000 &  1\,352   & 2\,341\,229   & 38\,696      & 1\,928
 \\ 0.0500 &  3\,718   & 117\,184\,420 & 308\,806     & 4\,870
 \\ 0.0250 &  11\,726  & & 2\,465\,102    & 14\,030
 \\ 0.0125 &  40\,144  &&& 44\,752
 \\ 0.0062 &  146\,926 &&& 156\,142
 \\ 0.0031 &  561\,912 &&& 580\,344
 \\
 \end{tabular}
\caption{The number of function evaluations  for numerical integration using moment--fitting (MF), Monte-Carlo (MC), sub-triangulation (ST), and local parametrization (LP) algorithms to handle cut elements.
    The table shows results with $m=5$.} \label{table2.1}
 \end{table}

\subsection{Unfitted FEM}\label{s_unfitted}

In the next series of experiments we solve the Poisson equation with Neumann's boundary condition in the unit disc domain defined implicitly
as $\Omega=\{\bx\in\mathbb{R}^2\,:\,\phi(\bx)<0\}$ with $\phi=|\bx|-1$.
We are interested in solving \eqref{Poiss}
with $\alpha=1$ and the right-hand side given in polar coordinates by
\[
f = (a^2+r^{-2} + 1) \sin(a \theta) \sin(\theta) - a/r \cos(a r) sin(\theta) +
                (c^2+1) \cos(c r) + c/r \sin(c r).
\]
The corresponding solution is $u(r, \theta) = \sin(a r) \sin(\theta) + \cos(c r) $.
In experiments we set  $a = 7 \pi/2$ and $c = 3 \pi$ .

\begin{figure}[ht]
  \centering
    \includegraphics[width=0.49 \textwidth]{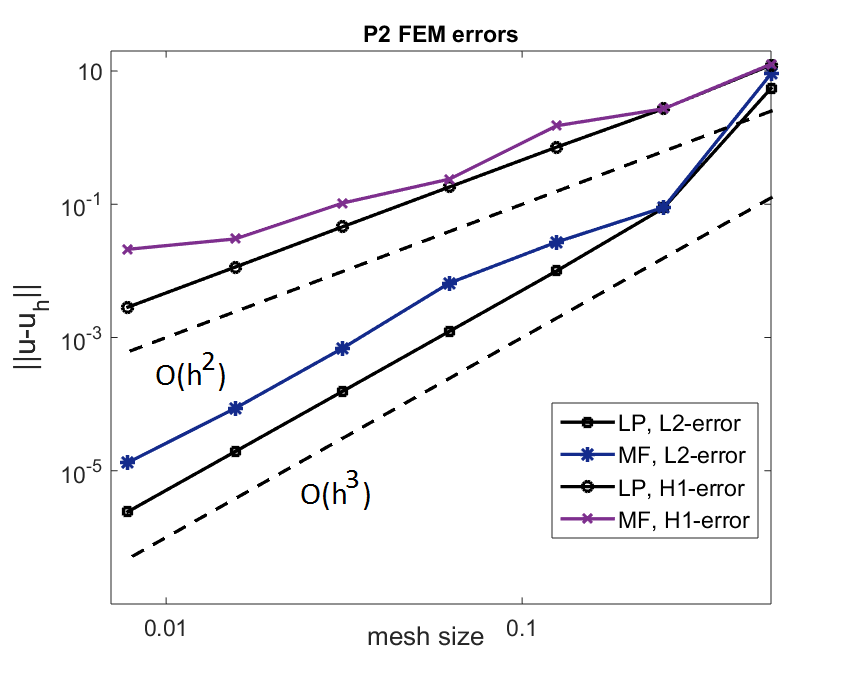}
    \includegraphics[width=0.49 \textwidth]{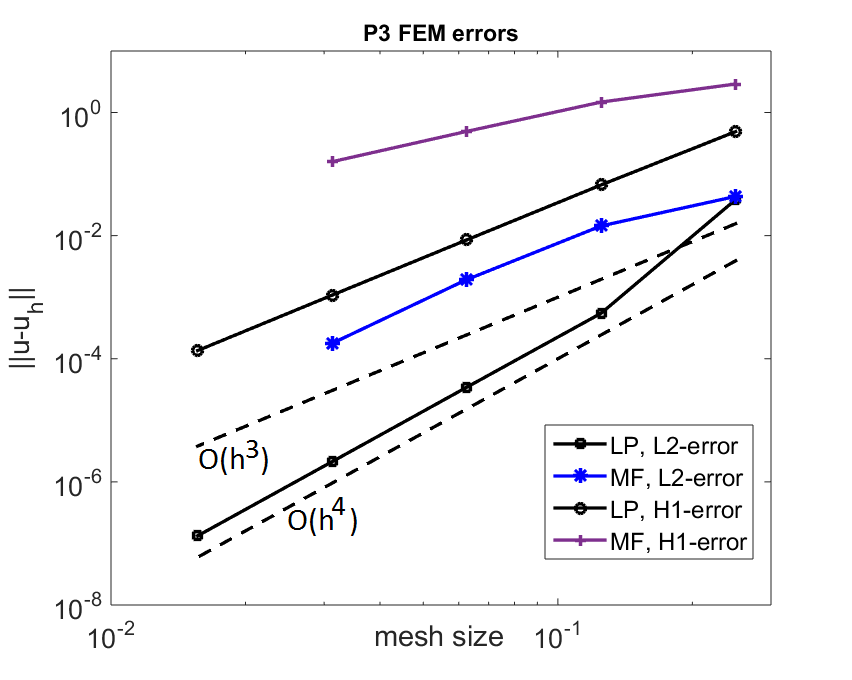}
    \caption{Finite element method error for $P_2$ (left) and $P_3$ (right) elements. The error plots are shown for moment--fitting (MF)
    and local parametrization (LP) used to treat cut elements. \label{fig_FEorder} }
\end{figure}

\begin{figure}[ht]
  \centering
    \includegraphics[width=0.49 \textwidth]{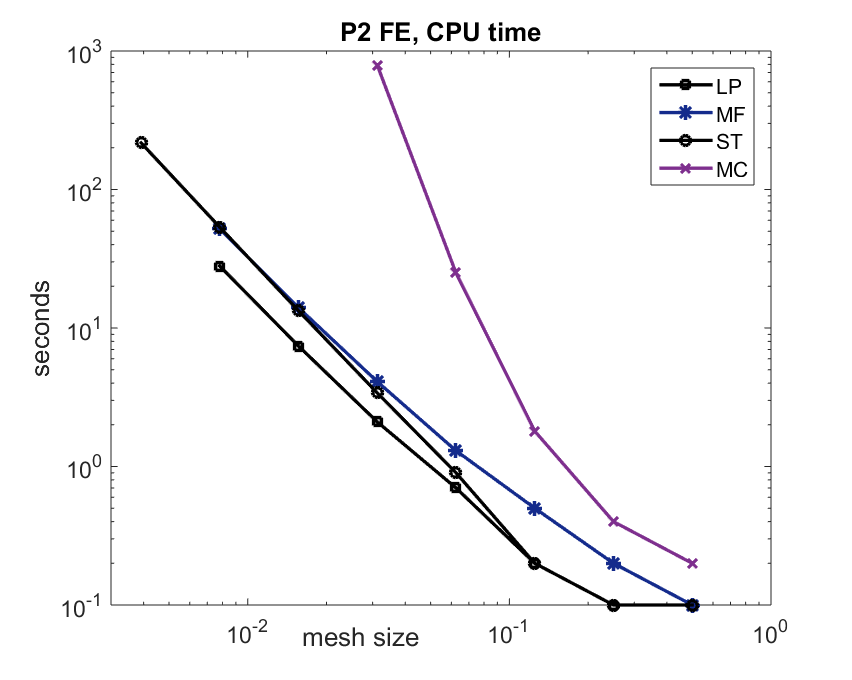}
    \includegraphics[width=0.49 \textwidth]{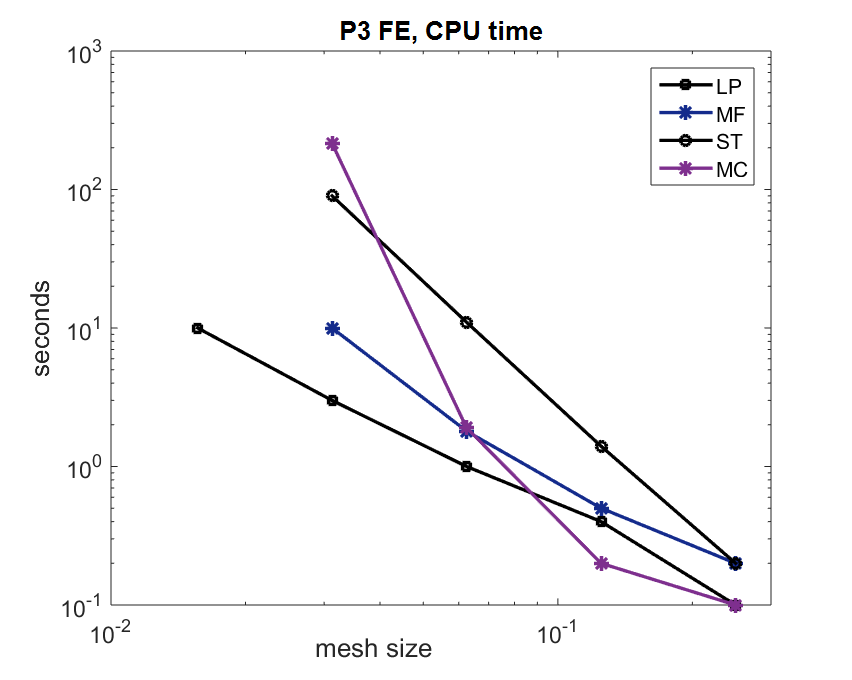}
    \caption{Dependence of the total CPU time to assemble stiffness matrices for $P_2$ (left) and $P_3$ (right) elements. Times are shown for moment--fitting (MF), local parametrization (LP), sub-triangulation (ST) and Monte-Carlo (MC) methods used to treat cut elements.\label{fig_FEtimes}}
\end{figure}

The bulk domain  $\Omega^{\rm bulk}=(-\frac32,\frac32)^2$ is triangulated using uniform meshes of sizes $h=0.5\times2^{-i}$, $i=0,\dots,8$.
We experiment with $P_2$ and $P_3$ finite elements. The discrete level set function is the finite element interpolant to the distance function $\phi$,
$\phi_h=J_h(\phi)$. Hence, the estimate \eqref{phi_h} holds with $q=2$ and $q=3$, respectively.
 To be consistent with the geometric error and polynomial order, all four integration methods for cut cells were set up to deliver the local error estimate \eqref{aux1} with $m=4$  and $m=5$, respectively. According to \eqref{FEerror1}, we should
 expect $O(h^2)$ and $O(h^3)$ convergence in the energy norm and $O(h^3)$ and $O(h^4)$ convergence in the $L^2(\Omega)$ norm.

Figure~\ref{fig_FEorder} shows the error plots for the unfitted finite element method \eqref{FEmeth} for different mesh sizes.
The results are shown with the moment--fitting and local parametrization algorithms used for the integration over cut triangles.
The FE errors for Monte-Carlo and sub-triangulation were very similar to those obtained  with the local parametrization quadratures and hence
they are not shown. These results are in perfect agreement with the error estimate \eqref{FEerror1}. Interesting that using the moment--fitting
method for computing the stiffness  matrix and the right-hand side leads to larger errors for the computed FE solution.

Further, Figure~\ref{fig_FEtimes} shows the CPU times needed for  the setup phase of the finite element method using different
numerical integration tools. For $P_2$ element, the moment--fitting, the local parametrization and the sub-triangulation method show similar
scaling with respect to $h$ since the complexity is dominated by the matrix assemble over internal triangles, while for $P_3$
elements local parametrization is superior in terms of final CPU times. As expected from the above analysis, the Monte--Carlo algorithm is non-optimal in either case.

\vskip.5cm

\subsection{Narrow-band unfitted FEM}
In the final series of experiments we apply the narrow-band unfitted finite element method \eqref{FEmethLB}
to solve the Laplace-Beltrami equation \eqref{LBeq} with $\alpha=1$ on the implicitly defined surface $\Gamma=\{\bx\in\mathbb{R}^2\,:\,\phi(\bx)=0\}$ with $\phi=|\bx|-1$.
The solution and  the right-hand side are given in polar coordinates by
\[
u = \cos(8 \theta),\quad f = 65 \cos(8 \theta).
\]

\begin{figure}[!ht]
  \centering
    \includegraphics[width=0.49 \textwidth]{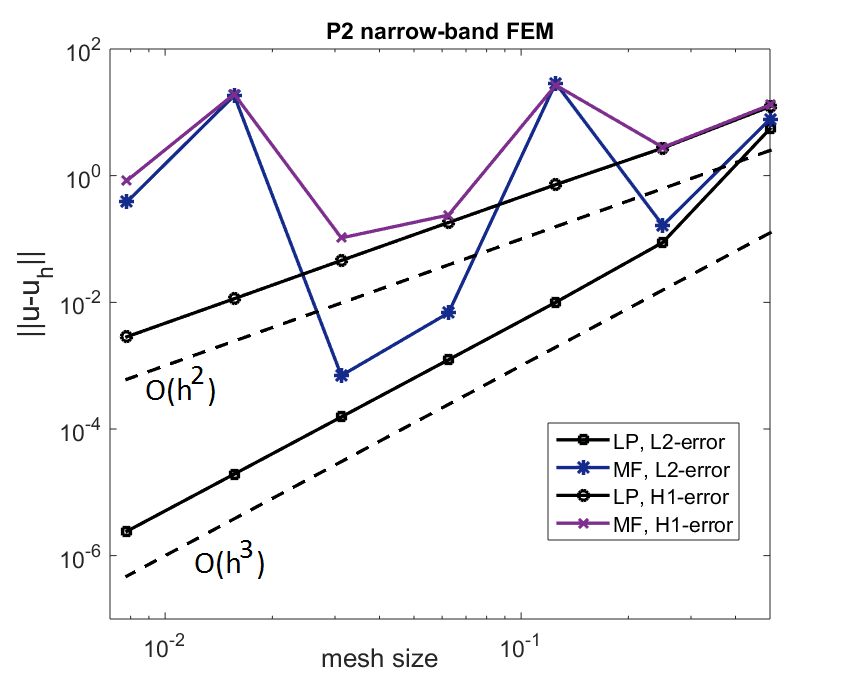}
    \includegraphics[width=0.49 \textwidth]{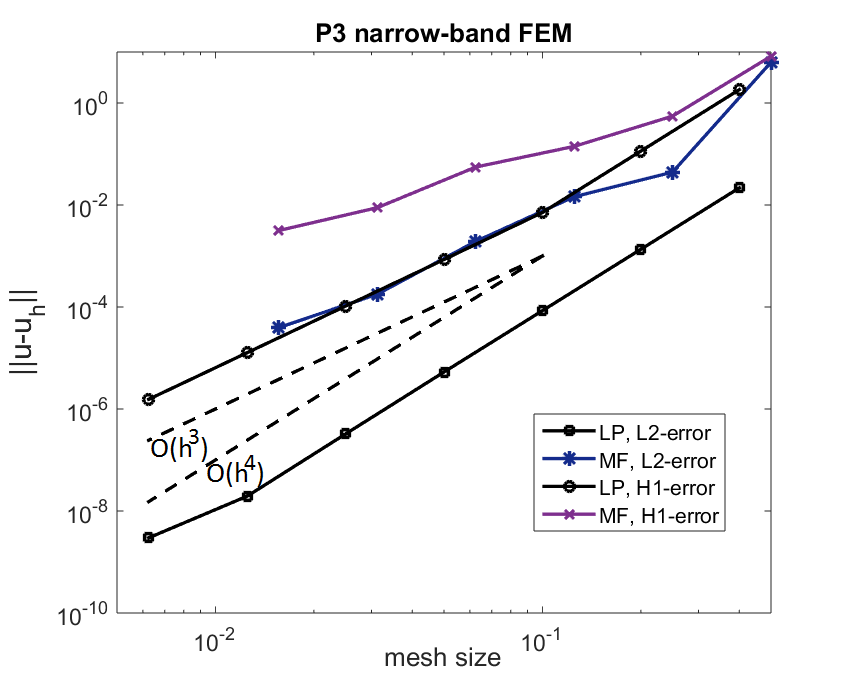}
    \caption{$L^2(\Gamma)$ and $H^1(\Gamma)$ errors for the narrow-band with $P_2$ (left) and $P_3$ (right) bulk elements. The error plots are shown for moment--fitting (MF)
    and local parametrization (LP) used to treat cut elements. \label{fig_FE2order} }
\end{figure}

The bulk domain,  triangulations and the discrete level set function are the same as used in the previous series of experiments in section~\ref{s_unfitted}.
For the extended finite element formulation \eqref{FEmethLB} we define the following narrow band domain:
\[
\Omega_h=\{\bx\in\mathbb{R}^2\,:\,|\phi_h(\bx)|<2h\}.
\]
The extension of the right-hand side is done along normal directions to $\Gamma$. For the discrete Hessian in  $\Omega_h$, we take the exact one
computed by $\bH_h:=\nabla^2\phi$ in $\Omega_h$.
Finite element space $V_h$ is the same as in section~\ref{s_unfitted} and is build on $P_2$ or $P_3$ piecewise polynomial continuous functions in $\Omega^{\rm bulk}$.
Similar to the previous test case,   four integration methods for cut cells were set up to deliver the local error estimate \eqref{aux1} with $m=4$  and $m=5$, respectively. According to \eqref{FEerror2}, we should
 expect $O(h^2)$ and $O(h^3)$ convergence in the energy norm.  Although there is no error estimate proved in the $L^2$ norm,
the optimal convergence order would be $O(h^3)$ and $O(h^4)$.

Figure~\ref{fig_FE2order} shows the error plots for the narrow-band  unfitted finite element method \eqref{FEmeth} for different mesh sizes.
All errors were computed over $\Gamma$ as stands in the estimate  \eqref{FEerror2}, rather than in the bulk.
The results are shown only for the moment--fitting and local parametrization algorithms used for the integration over cut triangles.
As before, the results with other methods were very similar to those obtained  with the local parametrization quadratures. The results
 with local parametrization are in perfect agreement with the error estimate \eqref{FEerror2} and predict the gain of one order in the $L^2(\Gamma)$ norm. Using the moment--fitting leads to unstable results in the case of quadratic finite elements and to sub-optimal convergence in the case of cubic elements. Note that in these two cases moment--fitting with quadratic and cubic basis $\mathcal{G}$, respectively, were used.

\section{Conclusions}\label{s_concl}
Building higher order quadrature rules for the numerical integration over implicitly defined curvilinear domains remains
a challenging problem, important in many applications of unfitted finite element methods. Well known approaches are not robust with respect to how
a surface cuts the mesh or have non-optimal computational complexity. In this paper we studied two methods of optimal complexity, namely,
the moment--fitting and the local parametrization. Although moment--fitting  delivers optimal accuracy for the integration of a smooth function over
a bulk curvilinear domain, its application to numerical PDEs were found to produce sub-optimal results. Local parametrization provides accurate and stable integration method. However, its extension to 3D problems is not straightforward and requires further studies. Developing more stable versions of the moment--fitting method, extending parametrization technique to higher dimensions, or devising ever different numerical approaches to \eqref{Int} all can be directions of further research.

\bibliographystyle{siam}
\bibliography{bibliography}
\end{document}